# ERGODIC BEHAVIOR OF LOCALLY REGULATED BRANCHING POPULATIONS

By M. Hutzenthaler and A. Wakolbinger

*Johann Wolfgang Goethe–Universität Frankfurt*

For a class of processes modeling the evolution of a spatially structured population with migration and a logistic local regulation of the reproduction dynamics, we show convergence to an upper invariant measure from a suitable class of initial distributions. It follows from recent work of Alison Etheridge that this upper invariant measure is nontrivial for sufficiently large super-criticality in the reproduction. For sufficiently small super-criticality, we prove local extinction by comparison with a mean field model. This latter result extends also to more general local reproduction regulations.

**1. Introduction.** In naturally reproducing populations one usually encounters an average number of more than one offspring per individual. However, given nonextinction, classical super-critical branching processes grow beyond all bounds. This is unrealistic because of bounded resources.

An efficient counteraction to unbounded population growth is achieved by a population-size dependent regulation of the reproduction dynamics. An example is the so-called *logistic branching process* [15] in which, in addition to the "natural" births and deaths in a super-critical branching mechanism, there are deaths resulting from a competition between any two individuals in the population. In Feller's diffusion limit, this leads to a negative drift term which is proportional to the squared population size and prevents the population size $X$ from escaping to $\infty$. This results in an eventual trapping of $X$ in zero.

An attempt to combat this extinction is to consider infinite populations modeled by a spatially extended version of the logistic branching process, with subpopulations living in discrete demes arranged in the $d$-dimensional lattice $\mathbb{Z}^d$, and with a (homogeneous) migration between the demes. This

---

leads to the following system of *interacting Feller diffusions with logistic growth* with local population sizes $X(i)$:

$$
\begin{aligned}
dX_t(i) = {} & \alpha\left(\sum_{j\in\mathbb{Z}^d} m(i,j)X_t(j) - X_t(i)\right)dt \\
& + \gamma X_t(i)(K - X_t(i))\,dt + \sqrt{2\beta X_t(i)}\,dB_t(i), \qquad i\in\mathbb{Z}^d.
\end{aligned}
\tag{1}
$$

Here, the $B(i)$ are independent standard Wiener processes, $m$ is the transition matrix of a random walk on $\mathbb{Z}^d$, and $\alpha,\beta,\gamma$ are nonnegative constants describing the rates of migration, branching and competition. The constant $K$ is called the *capacity*; it is the ratio of the *coefficient of super-criticality* $\gamma K$ and the competition rate $\gamma$.

Models with competition have been studied by various authors: Mueller and Tribe [18] and Horridge and Tribe [11] investigated an SPDE analogue of (1), with $d=1$ and $\mathbb{R}^1$ instead of $\mathbb{Z}^1$, and Etheridge [7], motivated by the work of Bolker and Pacala [2], investigated system (1) and its measure-valued analogue (with $\mathbb{Z}^d$ replaced by $\mathbb{R}^d$). These models also include long range competition. We emphasize that our methods make use of the fact that the interactions due to competition are solely within the same lattice site.

Using arguments from oriented percolation, Etheridge [7] shows that system (1) and also similar systems with nonlocal competition, when started from a spatially homogeneous initial state, do not suffer local extinction [i.e., $X_t(0)$ does not converge to zero as $t\to\infty$], provided the capacity $K$ is large enough. On the other hand, it was shown in the same paper by a coupling and comparison with subcritical branching (similar as in [18]) that a measure-valued analogue of (1) with certain nonlocal competition mechanisms suffers local extinction. The question whether lattice-based systems like (1) suffer local extinction for suitably small $K$ remained open. In the present paper we answer this question in the affirmative for the system (1) (Theorem 1). More precisely, we specify a strictly positive constant $\overline{K}$ such that for all capacities $K \le \overline{K}$ system (1) suffers local extinction. The constant $\overline{K}$ depends on the rates $\alpha,\beta$ and $\gamma$ of migration, branching and competition, but is uniform in all dimensions $d$ and migration matrices $m$.

The proof of local extinction is achieved by comparing (1) with a *mean field model* associated with (1). It turns out that there is no nontrivial equilibrium for (1) if 0 is the only equilibrium state of the mean field model. This model is mathematically very tractable; its critical capacity, which constitutes the upper bound $\overline{K}$, is the unique solution of the scalar equation appearing in Theorem 1.

We construct the maximal process $X^{(\infty)}$, which is the solution of (1) entering from infinity at time 0 (Theorem 2). As time tends to infinity, this process converges monotonically in distribution to the *upper invariant measure* of (1). Exploiting a self-duality of $X$ (Theorem 3), we show in Theorem 5 that the solution $X_t$ of (1), when started in a shift invariant nontrivial initial state, converges in distribution to this distinguished equilibrium.

Self-duality was used to prove ergodicity by other authors, for example, Horridge and Tribe [11] and Athreya and Swart [1]. In the latter paper, self-duality was established for the resampling selection model which is the solution of (1) where the Feller term $\sqrt{2\beta X_t(i)}$ is replaced by the Fisher–Wright term $\sqrt{2\beta X_t(i)(1 - X_t(i))}$ and where $K \leq 1$. Furthermore, Athreya and Swart study a branching coalescing particle model which in Feller's diffusion limit leads to the solution of (1). For both models, they prove existence of the maximal process and of the upper invariant measure.

**2. Main results.** We obtain the local extinction result and the result about existence of the maximal process and of the upper invariant measure for a more general class of *interacting locally regulated diffusions*, in which both the individual fertility and the branching rate depend on the local population size $X_t(i)$. The system of stochastic differential equations we are going to consider is

$$
\begin{aligned}
dX_t(i) = \alpha\left(\sum_{j \in G} m(i, j)X_t(j) - X_t(i)\right) dt \\
+ h(X_t(i))\, dt + \sqrt{2 \cdot g(X_t(i))}\, dB_t(i), \qquad i \in G,
\end{aligned}
\tag{2}
$$

where $G$ is an at most countable Abelian group. Unless stated otherwise, we will assume that the *migration matrix* $m$ is translation invariant, that is, $m(i, j) = m(0, j - i)$, that $\sum_{i \in G} m(0, i) = 1$, and that $m$ is irreducible, that is, $\forall i, j\ \exists n : m^{(n)}(i, j) > 0$. Assumptions on the functions $h$ and $g$ will be collected in Assumptions A1, A2 and A3 below.

An appropriate state space for (1) and (2) is provided by a construction going back to Liggett and Spitzer [17]: For given $m$, let $\sigma = (\sigma_i)_{i \in G}$ be summable and strictly positive such that

$$
\sum_i \sigma_i m(i, j) \leq C_{\mathrm{LS}}\sigma_j, \qquad j \in G,
\tag{3}
$$

for some $C_{\mathrm{LS}} < \infty$, and define the *Liggett–Spitzer space*

$$
\mathbb{E}_\sigma = \left\{\underline{x} \in \mathbb{R}_{\geq 0}^G : \|\underline{x}\|_\sigma := \sum_i \sigma_i |x_i| < \infty\right\}.
\tag{4}
$$

Notice that every translation invariant measure $\mu$ on $\mathbb{R}_{\geq 0}^G$ with $\int x_0 \mu(dx) < \infty$ is supported by $\mathbb{E}_\sigma$.

The following assumption guarantees existence and uniqueness of a strong $\mathbb{E}_\sigma$-valued solution of system (2).

ASSUMPTION A1. The functions $h: \mathbb{R}_{\geq 0} \to \mathbb{R}$ and $g: \mathbb{R}_{\geq 0} \to \mathbb{R}_{\geq 0}$ are locally Lipschitz continuous in $[0, \infty)$ and satisfy $h(0) = g(0) = 0$. In addition, the function $h$ is *upward Lipschitz continuous*, that is,

$$\text{sgn}(x - y)(h(x) - h(y)) \leq C|x - y| \tag{5}$$

for some constant $C$. Furthermore, $g$ is strictly positive on $(0, \infty)$ and satisfies the growth condition

$$\limsup_{x \to \infty} \frac{\sqrt{g(x)}}{x} < \infty. \tag{6}$$

PROPOSITION 2.1. *Assume Assumption* A1. *Then, for any* $\underline{x} \in \mathbb{E}_\sigma$, *the system* (2) *has a unique strong solution* $X = (X_t)$ *starting in* $\underline{x}$ *and having almost surely continuous paths in the normed space* $\mathbb{E}_\sigma$.

This will be proved in Section 3. Together with Assumption A1, the following condition on the drift function $h$ will be essential in Theorem 1.

ASSUMPTION A2. The function $h: \mathbb{R}_{\geq 0} \to \mathbb{R}$ is concave. Moreover, for some $x_0 > 0$, $h$ is negative on $[x_0, \infty)$, and satisfies

$$\int_{x_0}^{+\infty} \frac{1}{-h(x)} \, dx < \infty. \tag{7}$$

Note that (7) together with the upward Lipschitz continuity required in Assumption A1 implies

$$\lim_{x \to \infty} h(x) = -\infty. \tag{8}$$

For the interacting Feller diffusions with logistic growth (1), the functions $h$ and $g$ are of the form

$$h(x) = \gamma x(K - x), \qquad g(x) = \beta x. \tag{9}$$

In this case, Assumptions A1 and A2 are clearly satisfied.

THEOREM 1. *Assume Assumptions* A1 *and* A2. *Denote by* $X$ *the solution of equation* (2) *for an arbitrarily prescribed initial distribution on* $\mathbb{E}_\sigma$. *If*

$$\int_0^\infty \frac{h(y)}{g(y)} \exp\left( \int_1^y \frac{-\alpha x + h(x)}{g(x)} \, dx \right) dy \leq 0, \tag{10}$$

*then the process suffers local extinction, that is,*

$$(X_t) \Longrightarrow \delta_{\underline{0}} \qquad \text{as } t \to \infty. \tag{11}$$

*Here,* $\underline{0}$ *denotes the zero configuration.*

In the logistic Feller case (9), condition (10) simplifies to

$$(12) \qquad \int_0^\infty \exp\left(K\gamma y - \frac{\gamma\beta}{2}y^2\right) \cdot \alpha \exp(-\alpha y)\,dy \leq 1;$$

see the proof of Corollary 6 at the end of Section 5.

A first glance at system (1) might tempt one to believe that, even for small capacities $K$ (and $\alpha$ fixed), a suitably mobile migration $m$ in the dynamics (1) could prevent the system from suffering local extinction. However, Theorem 1 and condition (12) reveal that this is not the case.

The proof of Theorem 1 will be given in Section 6. Its main idea is a comparison with a *mean field model* corresponding to (2), given by the solution $V$ of

$$(13) \qquad dV_t = \alpha(\mathbf{E}V_t - V_t)\,dt + h(V_t)\,dt + \sqrt{2 \cdot g(V_t)}\,dB_t.$$

For this comparison, we will suppose that the initial condition $X_0$ of (2) has a distribution $\bar{\mu}$ which is *associated*, that is,

$$(14) \qquad \int f_1 \cdot f_2\,d\bar{\mu} \geq \int f_1\,d\bar{\mu} \int f_2\,d\bar{\mu}$$

for all bounded, coordinate-wise nondecreasing $f_1, f_2 : \mathbb{E}_\sigma \to \mathbb{R}$.

For the solution $X$ of (2), we will show that, for all $t \geq 0$, the marginal distributions of $X_t$ are bounded by the distribution of $V_t$ in the $\leq_{\mathrm{icv}}$ – order (where "icv" stands for "increasing, concave"; see [21] for this and related notions). More precisely, in Section 6 we will prove the following.

PROPOSITION 2.2. *Assume Assumption* A1 *and concavity of $h$. Let $X$ be a solution of* (2) *whose initial distribution $\bar{\mu}$ is associated. Assume that the $X_0(i)$, $i \in G$, are identically distributed and have finite expectation. Let $\overline{V}_t = (\overline{V}_t(i))_{i \in G}$ be a system of processes coupled through the initial state $\overline{V}_0(i) = X_0(i)$, $i \in G$, but following independent mean field dynamics, that is, every $\overline{V}_t(i)$ solves equation* (13) *with Brownian motion $B(i)$, where the $B(i)$, $i \in G$, are independent. Then*

$$(15) \qquad \mathbf{E}^{\bar{\mu}}f(X_t) \leq \mathbf{E}^{\bar{\mu}}f(\overline{V}_t)$$

*for all bounded, coordinate-wise nondecreasing and concave functions $f : \mathbb{E}_\sigma \to \mathbb{R}$ depending on only finitely many coordinates.*

The following theorem, whose proof will be given in Section 4, provides for the existence of a distinguished equilibrium state of (2), called the *upper invariant measure*. For proving Theorem 2, we will exploit the following assumption, which relaxes Assumption A2. Condition (16) ensures that the drift is "sufficiently negative" for large values of $X_t(i)$ so that the process "comes down from $\infty$."

ASSUMPTION A3. There exists a function $\hat{h} \geq h$ such that, for some $x_0 > 0$, $\hat{h}$ is negative and concave on $[x_0, \infty)$ and satisfies

$$\int_{x_0}^{+\infty} \frac{1}{-\hat{h}(x)} \, dx < \infty. \tag{16}$$

To prepare for Theorem 2, we need a bit of notation. If $\mu_1, \mu_2$ are probability measures on a partially ordered set $S$, then we say that $\mu_1$ is *stochastically smaller than or equal to* $\mu_2$, and we write $\mu_1 \leq \mu_2$, if there exists a random pair $(Y_1, Y_2)$ with marginal laws $(Y_i) = \mu_i$, $i = 1, 2$ and $Y_1 \leq Y_2$. We say that a sequence of probability measures $\mu_n$ *increases stochastically* to a probability measure $\mu_\infty$, denoted by $\mu_i \uparrow \mu_\infty$, if there exists a random sequence $(Y_i)$ which a.s. increases to $Y_\infty$ and has marginal distributions $(Y_i) = \mu_i$, $i = 1, 2, \ldots, \infty$.

THEOREM 2. *Assume Assumptions A1 and A3. There exists an $\mathbb{E}_\sigma$-valued process $(X_t^{(\infty)})_{t>0}$ with the following properties:*

(a) *For each $\varepsilon > 0$, $(X_t^{(\infty)})_{t \geq \varepsilon}$ is a solution of (2) starting at time $t = \varepsilon$.*

(b) *The first moment of $\|X_t^{(\infty)}\|_\sigma$ is finite for every $t > 0$.*

(b) *Let $\underline{x}^{(n)} = (x_i^{(n)})_{i \in G}, n = 1, 2, \ldots,$ be an increasing sequence in $\mathbb{E}_\sigma$ such that, for all $i \in G$,*

$$x_i^{(n)} \uparrow \infty \qquad \text{as } n \to \infty. \tag{17}$$

*If $X^{(n)}$ is the solution of (2) starting in $\underline{x}^{(n)} \in \mathbb{E}_\sigma$ at time zero, then*

$$\mathcal{L}(X_t^{(n)}) \uparrow \mathcal{L}(X_t^{(\infty)}) \qquad \text{as } n \uparrow \infty \ (t > 0). \tag{18}$$

(d) *There exists an equilibrium distribution $\bar{\nu}$ (called the* upper invariant measure*) for the dynamics (2) such that*

$$\mathcal{L}(X_t^{(\infty)}) \downarrow \bar{\nu} \qquad \text{as } t \uparrow \infty. \tag{19}$$

(e) *Any $\mathbb{E}_\sigma$-valued solution $X$ of (2) satisfies*

$$\mathcal{L}(X_t) \leq \mathcal{L}(X_t^{(\infty)}) \qquad (t > 0). \tag{20}$$

*In particular, any equilibrium $\nu$ is stochastically smaller than or equal to $\bar{\nu}$.*

(f) *Both the upper invariant measure $\bar{\nu}$ and $\mathcal{L}(X_t^{(\infty)})$ are translation invariant and associated.*

In the following two theorems, we exploit the specific form of the dynamics (1) of the interacting Feller diffusions with logistic growth. As it turns out, the solution of equation (1) has a property of *self-duality* which is helpful for the investigation of convergence to equilibria.

For the formulation of the self-duality result, write $m^\dagger$ for the transpose of the matrix $m$, choose a $\sigma^\dagger$ satisfying (3) with $m^\dagger$ instead of $m$, and recall that $\mathbb{E}_{\sigma^\dagger}$ denotes the corresponding Liggett–Spitzer space.

THEOREM 3. *Assume $\beta > 0$. Let $X$ and $X^\dagger$ be solutions of (1) with migration kernels $m$ and $m^\dagger$. Then we have the following self-duality:*

$$\mathbf{E}^{\underline{x}} \exp\left(-\frac{\gamma}{\beta} \langle X_t, \underline{y} \rangle\right) = \mathbf{E}^{\underline{y}} \exp\left(-\frac{\gamma}{\beta} \langle \underline{x}, X_t^\dagger \rangle\right) \tag{21}$$

*for all $\underline{x} \in \mathbb{E}_\sigma, \underline{y} \in \mathbb{E}_{\sigma^\dagger}$.*

A similar (though nonself-) duality for interacting Feller diffusions (or super- random walks), that is, (1) with $\gamma = 0$, is given by

$$\mathbf{E}^{\underline{x}} \exp(-\langle X_t, \underline{y} \rangle) = \exp(-\langle \underline{x}, v_t \rangle), \tag{22}$$

where $v = (v_t(i))$ solves the initial value problem

$$\frac{d}{dt} v_t(i) = \sum_{j \in G} m(i,j)(v_t(j) - v_t(i)) - v_t(i)^2, \qquad i \in G, v_0 = \underline{y}; \tag{23}$$

see, for example, Chapter 4 of [5].

The proof of Theorem 3 is contained in Section 7. The main advantage of the self-duality (21) is that instead of starting in a configuration with infinite total mass, we can analyze the evolution of the process started with finite total mass. For example, choose $y = \lambda \delta_0$ and $x$ with $x(i) \equiv$ const. Then the self-duality tells us that it makes no difference whether we study the law of $X_t(0)$ started in $x$, or that of the total mass $|X_t^\dagger| := \sum_i X_t^\dagger(i)$ with $X^\dagger$ started in $\lambda \delta_0$, $\lambda > 0$. This leads to the following corollary (see Lemma 7.1 together with Theorem 1):

COROLLARY 4. *Assume $\beta, \gamma > 0$. Let the parameters $\alpha, \beta, \gamma, K$ be such that inequality (12) holds. Then the solution $X$ of (1) started from an initial state of finite total mass [i.e., $\sum_i X_0(i) < \infty$] hits $\underline{0}$ in finite time a.s.*

In the more general situation of (2), self-duality is not available for proving the analogue of Corollary 4. We conjecture that also for interacting locally regulated diffusions (2), the local extinction (11) implies global extinction of finite mass processes as stated in Corollary 4, but we do not have a proof.

Theorem 3 will be the principal tool for proving convergence to the upper invariant measure specified in Theorem 2. This convergence will be the subject of Theorem 5 below. On an intuitive level, the reason for this convergence is as follows: There are two forces working against each other, super-critical branching and individual competition. The third ingredient is migration, which is important for spreading out newly produced mass. Super-critical branching increases mass, whereas competition among the individuals decreases it. If a (local) population size is large, then competition is stronger, whereas, as long as a local population size is small, the competition is negligible in comparison to the mass producing branching. Thus, there should be some attracting equilibrium state in which the two forces balance each other. This is the upper invariant measure.

THEOREM 5.  *Assume $\beta, \gamma > 0$. Let $X$ be a solution of (1) and suppose that $\mathcal{L}(X_0) \geq \mu$, where $\mu$ is a measure on $\mathbb{E}_\sigma$ which is translation invariant and does not charge the zero configuration $\underline{0}$. Then*

$$\mathcal{L}(X_t) \Longrightarrow \bar{\nu} \qquad as \ t \to \infty, \tag{24}$$

*where $\bar{\nu}$ is the upper invariant measure.*

From this it is clear that the only extremal translation invariant equilibrium distributions are $\delta_{\underline{0}}$ and $\bar{\nu}$. They coincide in case of local extinction and differ in case of survival. Section 8 will be devoted to the proof of Theorem 5.

**3. Preliminaries.** We turn first to the proof of existence and uniqueness of (2) under Assumption A1. Let us define

$$b_i(\underline{x}) := \alpha \left( \sum_{j \in G} m(i,j) x_j - x_i \right) + h(x_i), \qquad \underline{x} \in \mathbb{E}_\sigma, \tag{25}$$

where $\sigma = (\sigma_i)_{i \in G}$ satisfies (3). Denote $z^+ := z \vee 0$. By inequality (3) and Assumption A1, there exists a finite constant $C_1$ such that

$$\sum_{i \in M} \sigma_i \mathbb{1}_{x_i - y_i \geq 0} (b_i(\underline{x}) - b_i(\underline{y})) \leq C_1 \|(x. - y.)^+\|_\sigma \qquad \forall \underline{x}, \underline{y} \in \mathbb{E}_\sigma \tag{26}$$

for every subset $M \subseteq G$. From inequality (26), we will obtain monotonicity in the initial configuration. This monotonicity is a crucial property which we will exploit several times. First, we prove boundedness of second moments.

LEMMA 3.1.  *Suppose that $h$ and $g$ satisfy Assumption A1. Let $(X_t)$ be any weak solution of equation (2) with $\mathbf{E}\|X_0\|_\sigma^2 < \infty$ whose paths are continuous in $\mathbb{E}_\sigma$. Then there exists a constant $C < \infty$ such that, for each $T \geq 0$,*

$$\sup_{t \leq T} \mathbf{E}\|X_t\|_\sigma^2 \leq (1 + \mathbf{E}\|X_0\|_\sigma^2) e^{\mathrm{C}T} < \infty. \tag{27}$$

PROOF. Let $G_k$ be finite subsets of $G$ which monotonically exhaust $G$ as $k \to \infty$. Denote $\|x\|_{\sigma,k} := \sum_{i \in G_k} \sigma_i |x_i|$. Applying Itô's formula, we obtain

$$
(28) \quad
\begin{aligned}
d\|X_t\|_{\sigma,k}^2 &= 2\|X_t\|_{\sigma,k} \sum_{i \in G_k} \sigma_i (b_i(X_t)\,dt + \sqrt{2g(X_t(i))}\,dB_t(i)) \\
&\quad + 2 \sum_{i \in G_k} \sigma_i^2 g(X_t(i))\,dt.
\end{aligned}
$$

Let $n \in \mathbb{N}$. The continuous function $g$ is bounded on the interval $[0, n/\sigma_i]$ for every $i \in G_k$. Thus, the stochastic integrals on the right-hand side of (28) are $L_2$-martingales when stopped at time $\tau_n := \inf_{t \geq 0}\{\|X_t\|_\sigma \geq n\}$. By path continuity, we have $\tau_n \to \infty$ as $n \to \infty$ almost surely. Taking expectations, inequality (26) with $\underline{y} = \underline{0}$ implies

$$
(29) \quad
\begin{aligned}
&\mathbf{E}\|X_{t \wedge \tau_n}\|_{\sigma,k}^2 \\
&\leq \mathbf{E}\|X_0\|_{\sigma,k}^2 + 2\mathbf{E} \int_0^{t \wedge \tau_n} \left( C_1 \|X_s\|_\sigma^2 + \sum_{i \in G} \sigma_i^2 g(X_s(i)) \right) ds.
\end{aligned}
$$

By the growth condition (6), we know that $g(x) \leq C_2(1 + x^2)$ for some constant $C_2 < \infty$. Letting $k \to \infty$ and using monotone convergence, we obtain

$$
(30) \quad \mathbf{E}\|X_{t \wedge \tau_n}\|_\sigma^2 \leq \mathbf{E}\|X_0\|_\sigma^2 + C_3 \int_0^t (1 + \mathbf{E}\|X_{s \wedge \tau_n}\|_\sigma^2)\,ds
$$

for some constant $C_3 < \infty$. Applying Gronwall's inequality to the function $t \mapsto 1 + \mathbf{E}\|X_{t \wedge \tau_n}\|_\sigma^2$, we arrive at

$$
(31) \quad \mathbf{E}\|X_{t \wedge \tau_n}\|_\sigma^2 \leq (1 + \mathbf{E}\|X_0\|_\sigma^2)e^{C_3 t} - 1.
$$

Letting $n \to \infty$, Fatou's lemma completes the proof. $\square$

In the proof of Proposition 2.1, we need a stronger uniformity than Lemma 3.1 provides.

LEMMA 3.2. *Assume Assumption* A1. *Let* $(X_t)$ *be any weak solution of equation* (2) *satisfying condition* (27). *Then for each* $T \geq 0$, *there exists a constant* $\widetilde{C}_T < \infty$ *such that*

$$
(32) \quad \mathbf{E}\sup_{t \leq T} \|X_t\|_\sigma \leq \widetilde{C}_T(1 + \mathbf{E}\|X_0\|_\sigma + \mathbf{E}\|X_0\|_\sigma^2) < \infty.
$$

PROOF. Recall the definition of $G_k$ and $\|\cdot\|_{\sigma,k}$ from the proof of Lemma 3.1. Multiplying by $\sigma_i$ and summing over $i \in G_k$ in (2), we obtain, for $t \leq T$,

$$
(33) \quad
\begin{aligned}
&\|X_t\|_{\sigma,k} - \|X_0\|_{\sigma,k} \\
&= \int_0^t \sum_{i \in G_k} \sigma_i b_i(X_s)\,ds + \int_0^t \sum_{i \in G_k} \sigma_i \sqrt{2g(X_s(i))}\,dB_s(i).
\end{aligned}
$$

The estimate (26) implies that $\sum_{i \in G_k} \sigma_i b_i(X_s) \leq C_1 \|X_s\|_\sigma$. Thus, denoting the rightmost term in (33) by $M_t^k$, we obtain

$$(34) \qquad \sup_{u \leq t} \|X_u\|_{\sigma,k} \leq \|X_0\|_\sigma + \int_0^t C_1 \sup_{r \leq s} \|X_r\|_\sigma \, ds + \sup_{u \leq T} |M_u^k|.$$

The process $(M_t^k)$ is an $L_2$-martingale since, by the assumption $g(x) \leq C(1+x^2)$ and condition (27), the integrands $\sqrt{2g(X_s(i))}$ in (33) are square integrable, and the second moment $\mathbf{E}|M_T^k|^2 = \int_0^T 2 \sum_{i \in G_k} \sigma_i^2 \mathbf{E} g(X_s(i)) \, ds$ is bounded by $\overline{C}_T(1 + \mathbf{E}\|X_0\|_\sigma^2)$ for some constant $\overline{C}_T$. Thus, using the estimate $z \leq 1 + z^2$, we conclude from Doob's $L_2$-inequality that

$$(35) \qquad \mathbf{E} \sup_{u \leq T} |M_u^k| \leq 1 + \mathbf{E}|M_T^k|^2 \leq 1 + \overline{C}_T(1 + \mathbf{E}\|X_0\|_\sigma^2).$$

Therefore, taking expectations in (34) and applying monotone convergence, we obtain

$$(36) \qquad \begin{aligned} \mathbf{E} \sup_{u \leq t} \|X_u\|_\sigma &\leq \mathbf{E}\|X_0\|_\sigma + C_1 \int_0^t \mathbf{E} \sup_{r \leq s} \|X_r\|_\sigma \, ds \\ &\quad + 1 + \overline{C}_T(1 + \mathbf{E}\|X_0\|_\sigma^2) \end{aligned}$$

for all $t \leq T$. Now the assertion follows from Gronwall's inequality. $\quad\square$

The following monotone coupling lemma will be an important tool.

LEMMA 3.3. *Let $h_1, h_2$ and $g$ satisfy Assumption* A1, *and let $B = (B(i))_{i \in G}$ be a system of independent Brownian motions defined on some filtered probability space. For $\iota = 1, 2$, assume that $X^\iota$ is defined on the same probability space, satisfies equation* (2) *with Brownian motions $B(i)$, drift function $h_\iota$ and initial configuration $\underline{x}^\iota \in \mathbb{E}_\sigma$, and has continuous paths in $\mathbb{E}_\sigma$. Then*

$$(37) \qquad h_1 \leq h_2 \quad \text{together with} \quad \underline{x}^1 \leq \underline{x}^2 \text{ implies} \quad X_t^1 \leq X_t^2 \qquad \forall t \geq 0 \ a.s.$$

PROOF. The first part of the proof follows that of Theorem IV.3.2 in [12]. Let $1 > a_1 > \cdots > a_n > \cdots > 0$ be defined by

$$(38) \qquad \int_{a_1}^1 \frac{1}{u} \, du = 1, \int_{a_2}^{a_1} \frac{1}{u} \, du = 2, \ldots, \int_{a_n}^{a_{n-1}} \frac{1}{u} \, du = n, \ldots.$$

Notice that $a_n \to 0$ as $n \to \infty$. For every $n = 1, 2, \ldots$, define a continuous function $\psi_n(u)$ with support in $(a_n, a_{n-1})$ such that

$$(39) \qquad 0 \leq \psi_n(u) \leq \frac{2}{nu} \quad \text{and} \quad \int_{a_n}^{a_{n-1}} \psi_n(u) \, du = 1.$$

Furthermore, define

$$
\phi_n(x) := \mathbb{1}_{x>0} \int_0^x dy \int_0^y \psi_n(u)\, du, \qquad x \in \mathbb{R}. \tag{40}
$$

These functions satisfy $\phi_n \in C^2(\mathbb{R})$, $|\phi_n'(x)| \le 1$, $\phi_n''(x) = \mathbb{1}_{x>0}\psi_n(x)$, $\phi_n(x) \le x^+$ and $\phi_n(x) \to x^+$ as $n \to \infty$. Fix $i \in G$ and let $\tau_k := \inf\{t \ge 0 : X_t^1(i) \vee X_t^2(i) \ge k\}$. Write $\Delta_t^i := X_t^1(i) - X_t^2(i)$ and let $b_i^\iota$ be as in equation (25) with $h$ replaced by $h_\iota$, $\iota = 1, 2$. By Itô's formula,

$$
\begin{aligned}
\phi_n(\Delta_{t\wedge\tau_k}^i) &- \phi_n(\Delta_0^i) \\
&= \int_0^{t\wedge\tau_k} \phi_n'(\Delta_s^i)[\sqrt{2g(X_s^1(i))} - \sqrt{2g(X_s^2(i))}]\, dB_s(i) \\
&\quad + \int_0^{t\wedge\tau_k} \phi_n'(\Delta_s^i)[b_i^1(X_s^1) - b_i^2(X_s^2)]\, ds \\
&\quad + \tfrac{1}{2} \int_0^{t\wedge\tau_k} \phi_n''(\Delta_s^i)[\sqrt{2g(X_s^1(i))} - \sqrt{2g(X_s^2(i))}]^2\, ds.
\end{aligned} \tag{41}
$$

As $n \to \infty$, the left-hand side converges to $(\Delta_{t\wedge\tau_k}^i)^+ - (\Delta_0^i)^+$ in $L_1$ by dominated convergence and Lemma 3.1. In the rest of the proof, $C_1, C_2, \ldots$ will be suitably chosen finite constants. By Assumption A1, there exists a constant $C_1$ such that $g(x) \le C_1(1 + x^2)$. Thus, Lemma 3.1 implies $\mathbf{E}g(X_t^i(i)) < \infty$ and we have, by dominated convergence,

$$
\mathbf{E} \int_0^{t\wedge\tau_k} (\mathbb{1}_{\Delta_s^i>0} - \phi_n'(\Delta_s^i))^2 (\sqrt{2g(X_s^1(i))} - \sqrt{2g(X_s^2(i))})^2\, ds \xrightarrow{n\to\infty} 0. \tag{42}
$$

Hence, the first (stochastic) integral on the right-hand side converges in $L_2$ to the same expression with $\phi_n'(x)$ replaced by $\mathbb{1}_{x>0}$. For the second integral, notice that $b_i^\iota$ is globally Lipschitz continuous on $\{\underline{x} : x_i \le k\}$. Thus, for $s \le \tau_k$, $|b_i^\iota(X_s^\iota)|$ is bounded by $C_2\|X_s^\iota\|_\sigma$, which has finite expectation by Lemma 3.1, and we obtain, by dominated convergence,

$$
\int_0^{t\wedge\tau_k} |\mathbb{1}_{\Delta_s^i>0} - \phi_n'(\Delta_s^i)| \cdot |b_i^1(X_s^1) - b_i^2(X_s^2)|\, ds \xrightarrow{n\to\infty} 0. \tag{43}
$$

Finally, we consider the third integral on the right-hand side of equation (41). The local Lipschitz continuity of $g$ implies that $\sqrt{g}$ is globally $1/2$-Hölder continuous on the interval $[0, k]$. Therefore, the last integral in (41) is bounded by

$$
\int_0^{t\wedge\tau_k} \frac{2}{n|\Delta_s^i|} \cdot C_3|\Delta_s^i| \le \frac{2C_3 t}{n} \to 0 \qquad \text{as } n \to \infty. \tag{44}
$$

Putting these calculations together, equation (41) implies

$$
(\Delta_t^i)^+ - (\Delta_0^i)^+ = \int_0^t \mathbb{1}_{\Delta_s^i>0}[\sqrt{2g(X_s^1(i))} - \sqrt{2g(X_s^2(i))}]\, dB_s(i)
$$

$$(45) \qquad + \int_0^t \mathbb{1}_{\Delta_s^i > 0}[b_i^1(X_s^1) - b_i^2(X_s^2)]\,ds$$

for all $t \le \tau_k$ almost surely. By path continuity, we have $\tau_k \to \infty$ almost surely as $k \to \infty$ and thus, equation (45) holds for all $t \ge 0$. The stochastic integral on the right-hand side is an $L_2$-martingale because of $g(x) \le C_1(1 + x^2)$ and Lemma 3.1. Taking expectations, we arrive at

$$
\begin{aligned}
(46) \qquad &\mathbf{E}\Bigg[\sum_{i \in G} \sigma_i (\Delta_t^i)^+ - \sum_{i \in G} \sigma_i (\Delta_0^i)^+\Bigg] \\
&= \int_0^t \mathbf{E}\sum_{i \in G} \sigma_i \mathbb{1}_{\Delta_s^i > 0}[b_i^1(X_s^1) - b_i^2(X_s^2)]\,ds \\
&\le C_4 \int_0^t \mathbf{E}\sum_{i \in G} \sigma_i (\Delta_s^i)^+\,ds.
\end{aligned}
$$

In the last step, we used $b_i^1 \le b_i^2$ and inequality (26). By Gronwall's inequality, we obtain

$$(47) \qquad \mathbf{E}\|(\Delta_t^i)^+\|_\sigma \le \|(\Delta_0^i)^+\|_\sigma e^{C_4 t}, \qquad i \in G.$$

For later use, we note that this inequality implies

$$(48) \qquad \mathbf{E}\|\Delta_t^i\|_\sigma \le \|\Delta_0^i\|_\sigma e^{C_4 t}, \qquad i \in G,$$

if $b_1 = b_2$. For this, notice that $|x_i^1 - x_i^2| = (x_i^1 - x_i^2)^+ + (x_i^2 - x_i^1)^+$. The right-hand side of inequality (47) is zero by the assumption $x^1 \le x^2$, which finishes the proof of the monotonicity result for fixed $t \ge 0$. Finally, $X_t^1 \le X_t^2$ follows for all $t \in \mathbb{Q}_{\ge 0}$ and then by continuity of paths for all $t \ge 0$ almost surely. $\square$

PROOF OF PROPOSITION 2.1. Let $B = (B_i)_{i \in G}$ be a system of independent Brownian motions, and fix an initial condition $\underline{x} \in \mathbb{E}_\sigma$. We will prove existence of a solution of (2) similarly as in [9], where the system (2) is studied in the case $h = 0$. To this end, for finite $\Lambda \subseteq G$ and $i, j \in G$, we define $m^\Lambda(i, j) := m(i, j)\mathbb{1}_{i, j \in \Lambda}$ and consider the finite-dimensional system

$$
\begin{aligned}
(49) \qquad dX_t^\Lambda(i) &= \alpha \sum_{j \in \Lambda} m^\Lambda(i, j) X_t^\Lambda(j)\,dt - \alpha X_t^\Lambda(i)\,dt \\
&\quad + h(X_t^\Lambda(i))\,dt + \sqrt{2 \cdot g(X_t^\Lambda(i))}\,dB_t(i), \qquad i \in \Lambda.
\end{aligned}
$$

Under Assumption A1, equation (49) has a unique solution $X^\Lambda$ starting in $(x_i)_{i \in \Lambda}$. We extend $X^\Lambda$ to an infinite sequence (still denoted by the same symbol) by putting $X_t^\Lambda(i) := 0$ for $i \in G \setminus \Lambda$. Following the arguments in the proof of Theorem 1 in [9], one can show that there exists a process $X = (X_t(i))$ arising as the monotone limit

$$(50) \qquad X_t^\Lambda(i) \uparrow X_t(i) \qquad \text{as } \Lambda \uparrow G.$$

To show that $X$ has a.s. continuous paths in $\mathbb{E}_\sigma$, we first note that, for each finite $\Lambda \subseteq G$, the process $X^\Lambda$, being a finite-dimensional diffusion, has a.s. continuous paths and therefore satisfies

$$(51) \qquad \lim_{\delta \to 0} \mathbf{P}\left( \sup_{|t-s| \le \delta, s, t \le T} \|X_t^\Lambda - X_s^\Lambda\|_\sigma \ge \varepsilon \right) = 0$$

for all $\varepsilon > 0$ and $T > 0$.

For all finite $\Lambda \subseteq G$, the process $X^\Lambda$ satisfies the assumptions of Lemma 3.1, with $m(i,j)$ in (2) replaced by $m^\Lambda(i,j)$. Consequently, $X^\Lambda$ also satisfies (27), where the constant $C$ can be chosen uniformly in $\Lambda$. Therefore, by the monotone convergence (50), $X$ satisfies (27) and, due to Lemma 3.2, also (32).

Next, we set out to show that, for all $\varepsilon > 0$ and $T \ge 0$,

$$(52) \qquad \lim_{\Lambda \uparrow G} \mathbf{P}\left( \sup_{t \le T} \|X_t - X_t^\Lambda\|_\sigma \ge \varepsilon \right) = 0.$$

For this purpose, let $G_k$ and $\|x\|_{\sigma,k}$ be as in the proof of Lemma 3.1. From (50) together with the a.s. component-wise continuity of $X$ and Dini's theorem, we conclude that, for all $T > 0$ and $k \in \mathbb{N}$,

$$(53) \qquad \sup_{t \le T} \|X_t - X_t^\Lambda\|_{\sigma,k} \to 0 \qquad \text{a.s. as } \Lambda \uparrow G.$$

By (32) and dominated convergence, we therefore have

$$(54) \qquad \mathbf{E} \sup_{t \le T} \|X_t - X_t^\Lambda\|_{\sigma,k} \to 0 \qquad \text{a.s. as } \Lambda \uparrow G.$$

For every finite $\Lambda \subseteq G$ and $k \in \mathbb{N}$, we estimate

$$(55) \qquad \begin{aligned} &\mathbf{E} \sup_{t \le T} \|X_t - X_t^\Lambda\|_\sigma \\ &\qquad \le \mathbf{E} \sup_{t \le T} \|X_t - X_t^\Lambda\|_{\sigma,k} + 2\mathbf{E} \sup_{t \le T} \sum_{i \notin G_k} \sigma_i X_t(i). \end{aligned}$$

The rightmost term in (55) does not depend on $\Lambda$ and converges to 0, again because of (32) and dominated convergence. Together with (54) this implies that the left-hand side of (55) converges to zero, and proves (52).

For $\varepsilon, \delta$ and $T > 0$, we have the estimate

$$(56) \qquad \begin{aligned} &\mathbf{P}\left( \sup_{|t-s| \le \delta, s, t \le T} \|X_t - X_s\|_\sigma \ge 3\varepsilon \right) \\ &\qquad \le \mathbf{P}\left( \sup_{|t-s| \le \delta, s, t \le T} \|X_t^\Lambda - X_s^\Lambda\|_\sigma \ge \varepsilon \right) + 2\mathbf{P}\left( \sup_{t \le T} \|X_t - X_t^\Lambda\|_\sigma \ge \varepsilon \right). \end{aligned}$$

Because of (51) and (52) the left-hand side of (56) converges to 0 as $\delta \to 0$. This implies almost sure pathwise continuity.

For uniqueness, we proceed as follows. In the situation of Lemma 3.3, choose $h_1 = h_2$ and $x^1 = x^2$. Then pathwise uniqueness follows by applying Lemma 3.3 twice. Uniqueness in law and strong existence follow then from a Yamada–Watanabe type argument (see [22], Theorem 2.2). For the existence of a strong solution, it remains to show that the dependence of the unique solution on the initial configuration is measurable. This follows from the monotonicity result of Lemma 3.3.  □

LEMMA 3.4. *Let $h$ and $g$ satisfy Assumption* A1. *The strong solution $X_t$ of system* (2) *is monotonically continuous in its initial configuration in the following sense: Let $\underline{x}^{(n)}, \underline{x} \in \mathbb{E}_\sigma$, be the starting points of $X_t^{(n)}$ and $X_t$, such that*

$$(57) \qquad \underline{x}^{(n)} \uparrow (\downarrow) \underline{x} \qquad as \ n \uparrow \infty.$$

*Then*

$$(58) \qquad X_t^{(n)} \uparrow (\downarrow) X_t \qquad \forall t \geq 0 \ as \ n \uparrow \infty \ a.s.$$

PROOF. In equation (48), let $h^1 = h^2 := h$, $X_t^1 := X_t$ and $X_t^2 := X_t^{(n)}$. Letting $n \to \infty$, this implies $L_1$-convergence of $X_t - X_t^{(n)}$ for fixed time $t \geq 0$. The monotonicity result of Lemma 3.3 finishes the proof.  □

## 4. The upper invariant measure. Proof of Theorem 2.

PROOF OF THEOREM 2. To fix notation, let us write $\mathcal{L}^{\underline{x}}(X_t)$ for the distribution of $X_t$ [the solution of (2)] starting from an element $\underline{x} \in \mathbb{E}_\sigma$. For $N \in \mathbb{N}$, we define the element $\underline{N} \in \mathbb{E}_\sigma$ by $\underline{N}(i) \equiv N$, $i \in G$. Let $X_t^N$ be the process started from $\underline{N}$. By Lemma 3.3, the sequence $X_t^N$ is nondecreasing in $N$ for all $t > 0$; let us write $X_t^{(\infty)}$ for its a.s. limit.

Now let $(\underline{x}^{(n)})$ be a sequence as in Theorem 2(c). For all $n \in \mathbb{N}$, we conclude from Lemma 3.4 that

$$(59) \qquad \mathcal{L}(X_t^{(\infty)}) \underset{N \to \infty}{\nwarrow} \mathcal{L}^{\underline{N}}(X_t) \geq \mathcal{L}^{\underline{x}^{(n)} \wedge N}(X_t) \underset{N \to \infty}{\nearrow} \mathcal{L}^{\underline{x}^{(n)}}(X_t).$$

Again by Lemma 3.4, we obtain, for all $N \in \mathbb{N}$,

$$(60) \qquad \mathcal{L}^{\underline{x}^{(n)} \wedge N}(X_t) \underset{n \to \infty}{\nearrow} \mathcal{L}^{\underline{N}}(X_t).$$

Thus, by a diagonal argument, there is a subsequence $\underline{x}^{(n_N)}$ of $\underline{x}^{(n)}$ such that

$$(61) \qquad \mathcal{L}^{\underline{x}^{(n_N)} \wedge N}(X_t) \underset{N \to \infty}{\nearrow} \mathcal{L}(X_t^{(\infty)}).$$

Together with inequality (59) and monotonicity (Lemma 3.3), this results in

$$\mathcal{L}^{x^{(n_N)}}(X_t) \underset{N \to \infty}{\nearrow} \mathcal{L}(X_t^{(\infty)}).$$ (62)

As $(x^{(n)})$ is an increasing sequence, (62) is equivalent to (18).

The next step shows that the limit is finite almost surely. Let $\hat{h} \geq h$ be the function given by Assumption A3. Notice that $\hat{h}$ may be replaced by $\hat{h} + C$ for every constant $C \geq 0$. Furthermore, $h$ is bounded above. Thus, we may assume that $\hat{h} \geq h$ is concave. By Itô's formula, Lemma 3.1 and translation invariance,

$$\frac{d}{dt} \mathbf{E} X_t^N(i) = \mathbf{E} h(X_t^N(i)) \leq \mathbf{E} \hat{h}(X_t^N(i)) \leq \hat{h}(\mathbf{E} X_t^N(i)).$$ (63)

For the last step, we applied Jensen's inequality. Therefore, the expectation is bounded above by the deterministic function $y(t, x)$ satisfying

$$\frac{d}{dt} y(t, x) = \hat{h}(y(t, x)), \qquad y(0, x) = x.$$ (64)

The concave function $\hat{h}(x)$ converges to $-\infty$ as $x \to \infty$. Choose $x_0$ such that $\hat{h}$ is strictly negative for all $x \geq x_0$. Then for all $x > x_0$ and $t > 0$, we have $x_0 < y(t, x) < x$. From (64), we obtain by separation of variables that the solution satisfies

$$
\begin{aligned}
t &= -\int_{y(t,x)}^{x} \frac{1}{\hat{h}(z)} \, dz \\
&\leq \int_{y(t,x)}^{\infty} \frac{1}{-\hat{h}(z)} \, dz \downarrow \int_{\lim_{x \to \infty} y(t,x)}^{\infty} \frac{1}{-\hat{h}(z)} \, dz \qquad \text{as } x \to \infty.
\end{aligned}
$$ (65)

For the monotone convergence, notice that $y(t, x)$ is nondecreasing in $x$ and that all integrals are finite by inequality (16). Hence, if $\lim_{x \to \infty} y(t, x)$ was infinite for $t > 0$, then we would face the contradiction $0 < t \leq 0$. Therefore, we arrive at

$$
\begin{aligned}
\mathbf{E} \|X_t^{(\infty)}\|_\sigma &= \sum_{i \in G} \sigma_i \uparrow \lim_{N \to \infty} \mathbf{E} X_t^N(i) \\
&\leq \sum_{i \in G} \sigma_i \lim_{x \to \infty} y(t, x) < \infty, \qquad t > 0.
\end{aligned}
$$ (66)

From Lemma 3.4, it is then clear that, for all $\varepsilon > 0$, the solution of (2) which starts at time $t = \varepsilon$ from $X_\varepsilon^{(\infty)}$ is the a.s. monotone limit (as $N \to \infty$) of the solutions of (2) starting from $X_\varepsilon^{(N)}$ at time $\varepsilon$, or, equivalently, starting from $\underline{N}$ at time 0. At the beginning of the proof we defined $X_t^{(\infty)}$ as this limit; hence, we have so far proved parts (a), (b) and (c) of Theorem 2.

A similar argument as in (59) proves that the process with initial measure $\mu$ is dominated by the maximal process, which is part (e).

To prove part (d), fix $0 < s < t$. By part (e),

$$(67) \qquad \mathcal{L}(X_r^{(\infty)}) \geq \mathcal{L}^{X_{t-s}^{(\infty)}}(X_r).$$

Using this with $r = s$, we get the inequality

$$(68) \qquad \mathcal{L}(X_s^{(\infty)}) \geq \mathcal{L}^{X_{t-s}^{(\infty)}}(X_s) = \mathcal{L}(X_t^{(\infty)}),$$

where the last equality follows from the Markov property. We conclude from this monotonicity that $\mathcal{L}(X_t^{(\infty)}) \downarrow \bar{\nu}$ for some probability measure $\bar{\nu}$ on $\mathbb{E}_\sigma$, which by continuity in the initial configuration (Lemma 3.4) is an equilibrium distribution of the dynamics (2).

Next, we show that the upper invariant measure is translation invariant and is associated. Both properties are preserved under weak limits. Furthermore, we will argue that these properties are preserved under the dynamics. The constant configuration $X_0^N \equiv N$ is both translation invariant and associated. Hence, both $X_t^{\frac{N}{}}$ and $X_t^{(\infty)}$ have these properties for all $t > 0$. Therefore, the claim follows.

The translation invariance of the migration kernel implies that the dynamics (2) preserve translation invariance. To prove the preservation of associated measures, we will argue in a similar way as in [4] where the analogue of (2) with $h = 0$ and $[0,1]^G$ instead of $\mathbb{R}_{\geq 0}^G$ was treated. We first consider the approximation scheme $(X^\Lambda, \Lambda)$ with finite $\Lambda \subset G$, used to prove the existence part of Proposition 2.1. For fixed $\Lambda$, Theorem 1.1 in [10], together with a uniform approximation of $h$ and $g$ on compact intervals by smooth and bounded functions $h_k$ and $g_k$ with $\inf_{x \geq 0} g_k(x) > 0$, shows that, for an associated initial distribution $\mathcal{L}(X_0)$, the projections of $\mathcal{L}(X_t^\Lambda)$ to $\mathbb{R}_{\geq 0}^\Lambda$ are associated. Since $\mathcal{L}(X_t^\Lambda)$ approximates $\mathcal{L}(X_t)$ as $\Lambda \uparrow G$, the claim follows. □

## 5. The mean field model.
In this section we study the dynamics

$$(69) \qquad dV_t = \alpha(\mathbf{E}V_t - V_t)\,dt + h(V_t)\,dt + \sqrt{2g(V_t)}\,dB_t.$$

It can be shown (but will not be required for the subsequent proofs) that (69) arises as the limit of a sequence of processes following the dynamics (2), where $G$ is replaced by a finite set $G_n$ of cardinality $n$ and $m^{(n)}(i,j) = 1/n$ for $i, j \in G_n$. This type of limit is known as *mean field* or *Vlasov–McKean limit*; we will therefore address (69) briefly as *mean field model*. Intuitively, a uniform migration which spreads out mass as far as possible should be good for survival, and conversely, extinction of $V$ governed by (69) should imply extinction of $X$ governed by (2). With this motivation in mind, we investigate in this section conditions on $h$ and $g$ under which the dynamics (69) admit a nontrivial equilibrium distribution.

To this end, we consider the following:

LEMMA 5.1. *Suppose that Assumption* A1 *holds and that*

$$\text{(70)} \qquad \exists y_0 > 0 : h|_{[0,y_0]} \geq 0 \quad and \quad 0 \not\equiv h|_{[y_0,\infty)} \leq 0.$$

*There is no nontrivial invariant measure for the dynamics* (69) *if and only if*

$$\text{(71)} \qquad \int_0^\infty \frac{h(y)}{g(y)} \exp\left(\int_{y_0}^y \frac{-\alpha x + h(x)}{g(x)}\, dx\right) dy \leq 0.$$

*If condition* (71) *is not satisfied, then there is exactly one nontrivial invariant measure. Assume additionally that* $\lim_{\varepsilon \to 0} \int_\varepsilon^{y_0} \frac{-\alpha x + h(x)}{g(x)}\, dx$ *exists in* $(-\infty, \infty]$. *Then condition* (71) *is equivalent to*

$$\text{(72)} \qquad \int_0^\infty \frac{\alpha y}{g(y)} \exp\left(\int_0^y \frac{-\alpha x + h(x)}{g(x)}\, dx\right) dy \leq 1.$$

PROOF. Let $\theta > 0$ and consider the process given by

$$\text{(73)} \qquad dV_t^\theta = \alpha(\theta - V_t^\theta)\, dt + h(V_t^\theta)\, dt + \sqrt{2g(V_t^\theta)}\, dB_t.$$

By standard theory (e.g., pages 220f and 241 in [14]), the equilibrium distribution of (73) is

$$\text{(74)} \qquad \Gamma_\theta(dy) = \frac{C_\theta}{g(y)} \exp\left(\int_{y_0}^y \frac{\alpha(\theta - x) + h(x)}{g(x)}\, dx\right) dy =: C_\theta \Phi(y)\, dy,$$

where $C_\theta \in (0, \infty)$ is the normalizing constant. Indeed, existence of an equilibrium of (73) is clear since the drift in zero is positive in zero and becomes sufficiently negative near $\infty$; formally, this follows from the finiteness of the integral $\int_0^\infty \Phi(y)\, dy$, which can be checked easily.

Obviously, (69) admits a nontrivial equilibrium if and only if $\int y\Gamma_\theta(dy) = \theta$ has a positive solution. Hence, all we need to do is to characterize the situations where

$$\text{(75)} \qquad \not\exists \theta > 0 : f(\theta) := \alpha \int \frac{y - \theta}{C_\theta} \Gamma_\theta(dy) = 0.$$

We eliminate one occurrence of $\theta$ on the left-hand side of (75) by an integration by parts:

$$f(\theta) = \int_0^\infty \frac{\alpha(y - \theta)}{g(y)} \exp\left(\int_{y_0}^y \frac{\alpha(\theta - x)}{g(x)}\, dx\right) \exp\left(\int_{y_0}^y \frac{h(x)}{g(x)}\, dx\right) dy$$

$$(76) \qquad = \lim_{\varepsilon \to 0} \left[ \exp\left( \int_{y_0}^{y} \frac{\alpha(\theta - x) + h(x)}{g(x)} \, dx \right) \right]_{1/\varepsilon}^{\varepsilon}$$

$$+ \int_0^\infty \frac{h(y)}{g(y)} \left( \exp\left( \int_{y_0}^{y} \frac{\alpha}{g(x)} \, dx \right) \right)^{\theta} \exp\left( \int_{y_0}^{y} \frac{-\alpha x + h(x)}{g(x)} \, dx \right) dy.$$

We now analyze the two boundary terms on the right-hand side of (76). In the following calculations, $C_i$ are finite constants. Recall that $h$ is nonpositive for large arguments. Furthermore, in Assumption A1 we assumed $g(x) \leq Cx^2$ for some constant $C$ and all $x \geq y_0 > 0$. With this, the expression coming from the boundary value $1/\varepsilon$ tends to zero as $\varepsilon \to 0$:

$$0 \leq \exp\left( \int_{y_0}^{1/\varepsilon} \frac{\alpha(\theta - x) + h(x)}{g(x)} \, dx \right)$$

$$(77)$$

$$\leq C_1 \exp\left( \int_{y_0 \vee \theta}^{1/\varepsilon} \frac{\alpha(\theta - x)}{Cx^2} \, dx \right) \xrightarrow{\varepsilon \to 0} 0.$$

For the other boundary term, we recall that $h$ is nonnegative for small arguments and estimate

$$0 \leq \exp\left( - \int_{\varepsilon}^{y_0} \frac{\alpha(\theta - x) + h(x)}{g(x)} \, dx \right)$$

$$(78)$$

$$\leq C_2 \exp\left( - \int_{\varepsilon}^{y_0 \wedge (\theta/2)} \frac{\alpha \theta / 2}{g(x)} \, dx \right)$$

$$\xrightarrow{\varepsilon \to 0} C_3 \exp\left( -\theta \frac{\alpha}{2} \int_{0+} \frac{1}{g(x)} \, dx \right).$$

By assumption, $g$ is locally Lipschitz in zero and thus, $g(x) \leq C_4 x$ in a neighborhood of zero. Together with $\theta > 0$, this implies that all boundary terms vanish. Notice that the expression coming from the boundary value $\varepsilon$ does not need to be zero in case $\theta = 0$.

At this point we have seen that $f$ can be rewritten as

$$(79) \quad f(\theta) = \int_0^\infty \frac{h(y)}{g(y)} \left( \exp\left( \int_{y_0}^{y} \frac{\alpha}{g(x)} \, dx \right) \right)^{\theta} \exp\left( \int_{y_0}^{y} \frac{-\alpha x + h(x)}{g(x)} \, dx \right) dy$$

for $\theta > 0$. We will show that $f$ is strictly decreasing and continuous in $\theta > 0$. For this, consider the function

$$(80) \qquad\qquad \theta \mapsto h(y) \left( \exp\left( \int_{y_0}^{y} \frac{\alpha}{g(x)} \, dx \right) \right)^{\theta}$$

for fixed $y \geq 0$. If $y < y_0$, then $h(y) \geq 0$ and the integral is negative. If $y > y_0$, then $h(y) \leq 0$ and the integral is positive. In both situations, the function in (80) is nonincreasing. Furthermore, there is an interval $[y_1, y_2]$ with $y_0 \leq y_1 < y_2$ where $h(y) < 0$ and where the function in (80) is strictly decreasing and converging to $-\infty$. The integral over $[0, y_0]$ on the right-hand side of (79) is continuous and nonincreasing in $\theta > 0$, and bounded in $\theta \geq 1$. This follows from dominated convergence and the fact that the integral over $[0, \varepsilon]$ on the right-hand side of (79) is bounded above by

$$(81) \quad \int_0^\varepsilon \frac{1}{g(y)} \exp\left(\int_\varepsilon^y \frac{\alpha\bar{\theta}}{2g(x)}\, dx\right) dy \sup_{x \leq y_0} h(x) \exp\left(\int_{y_0}^\varepsilon \frac{\alpha(\bar{\theta}-x)+h(x)}{g(x)}\, dx\right)$$

$$\leq \exp\left(\int_\varepsilon^y \frac{\alpha\bar{\theta}}{2g(x)}\, dx\right)\Big|_0^\varepsilon \cdot C < \infty$$

for all $\theta \geq \bar{\theta} > 0$, where $\varepsilon > 0$ is such that $|\alpha x - h(x)| \leq \alpha\bar{\theta}/2$ for all $x \leq \varepsilon$. By monotone convergence, the integral over $[y_0, \infty)$ on the right-hand side of (79) is continuous and strictly decreasing in $\theta > 0$, and decreases to $-\infty$. Thus, the function $f$ is continuous and strictly decreasing in $\theta > 0$ with $f(\infty) = -\infty$. Hence, condition (75) is satisfied if and only if $\lim_{\theta \to 0} f(\theta) \leq 0$. Note that, by strict monotonicity of $f$, there is at most one nontrivial invariant measure.

For the limit $\theta \to 0$ in equation (79), we use monotone convergence (for the $\int_0^{y_0}$ part) and dominated convergence (for the $\int_{y_0}^\infty$ part). Thus, we have

$$(82) \quad \lim_{\theta \to 0} f(\theta) = \int_0^\infty \frac{h(y)}{g(y)} \exp\left(\int_{y_0}^y \frac{-\alpha x + h(x)}{g(x)}\, dx\right) dy.$$

Therefore, $\lim_{\theta \to 0} f(\theta) \leq 0$ is equivalent to condition (71).

Now, additionally assume that $\lim_{\varepsilon \to 0} \int_\varepsilon^{y_0} \frac{-\alpha x + h(x)}{g(x)}\, dx$ exists in $(-\infty, \infty]$. Then reversing the calculation in (76) with $\theta = 0$, we arrive at

$$(83) \quad \lim_{\theta \to 0} f(\theta) = \int_0^\infty \frac{\alpha y}{g(y)} \exp\left(\int_{y_0}^y \frac{-\alpha x + h(x)}{g(x)}\, dx\right) dy$$

$$- \exp\left(-\lim_{\varepsilon \to 0} \int_\varepsilon^{y_0} \frac{-\alpha x + h(x)}{g(x)}\, dx\right).$$

If the limit on the right-hand side is $\infty$, then $\lim_{\theta \to 0} f(\theta) > 0$ and a nontrivial invariant measure exists. The assertion is true in this case because the left-hand side of (72) is $\infty$. Otherwise, the limit on the right-hand side of (83) is finite. Then multiply the equation with $\exp\left(\int_0^{y_0} \frac{-\alpha x + h(x)}{g(x)}\, dx\right)$ and merge the two integrals $\int_0^{y_0}$ and $\int_{y_0}^y$ into one integral. Hence, we see that (71) and (72) are equivalent. □

We now specialize this result to the logistic Feller case, where condition (72) can be simplified.

COROLLARY 6. *Consider the mean field model* (69) *with* $h(x) = \gamma x(K - x)$ *and* $g(x) = \beta x$. *Assume* $\alpha, \gamma, \beta > 0$ *and let* $\overline{K} > 0$ *be uniquely determined by*

$$(84) \qquad \int_0^\infty \exp\left(\overline{K}\gamma y - \frac{\gamma\beta}{2}y^2\right) \cdot \alpha \exp(-\alpha y)\, dy = 1.$$

*There is no nontrivial invariant measure for* (69) *if and only if* $0 \le K \le \overline{K}$.

PROOF. First of all, convince yourself that Assumptions A1 and A2 hold. Thus, Lemma 5.1 applies if $K > 0$. After an integration and a change of variables $(y \to \beta y)$, condition (72) takes the form (12). The left-hand side in (12) is strictly increasing in $K$, tends to $\infty$, is continuous in $K$ by monotone convergence and is smaller than one for $K = 0$. Hence, $\overline{K}$ exists and is unique. By monotonicity, condition (12) holds if and only if $K \le \overline{K}$. □

For example, in the case $\alpha = \gamma = \beta = 1$ formula (84) gives the numerical value $\overline{K} = 0.6973\cdots$.

The following extinction result for the mean field dynamics is a fairly direct consequence of Lemma 5.1.

LEMMA 5.2. *Consider the mean field model given by* (69). *Suppose that Assumptions* A1, A3 *and condition* (70) *hold. Then inequality* (71) *implies local extinction:*

$$(85) \qquad \mathcal{L}(V_t) \Longrightarrow \delta_0 \qquad (as\ t \to \infty)$$

*for any initial law.*

PROOF. Paralleling the arguments in Section 4, one infers the existence of the maximal process $V^{(\infty)}$ for the dynamics (69), which obeys $\mathcal{L}(V_t^{(\infty)}) \ge \mathcal{L}(V_t)$. Again, this maximal process converges to an invariant measure. However, by Lemma 5.1 and condition (71), the trivial measure $\delta_0$ is the only invariant measure. This implies the assertion. □

**6. Comparison with the mean field model. Proof of Theorem 1.** The main idea for the proof of Theorem 1 is the assertion that the interacting locally regulated diffusions are dominated by the mean field model. The intuition behind this is that a uniform spread of mass reduces competition and therefore is good for survival, and that the mean field model arises as a limit of uniform migration models (see Section 5).

We proceed in two steps to prove Theorem 1. First, we establish a comparison between the system of interacting locally regulated diffusions (2) and the mean field model (69) which implies that it is more likely for the latter to survive. Then we exploit the fact (proved in Section 5) that for some parameter configurations not even the mean field model survives.

The proof of the comparison result will first treat the case where the functions $h$ and $g$ satisfy the following assumptions.

Assumption A4. *The set* $I$ *is a closed finite interval of the form* $[0, c]$, $0 < c < \infty$. *The functions* $h \colon I \to \mathbb{R}$ *and* $\sqrt{g} \colon I \to \mathbb{R}$ *are twice continuously differentiable on* $I$ *and satisfy* $h(0) = g(0) = g(c) = 0 > h(c)$. *Furthermore,* $g$ *is strictly positive on* $(0, c)$.

The proof of Proposition 2.2 is based on the following lemma.

Lemma 6.1. *Let* $h$ *and* $g$ *satisfy Assumption* A4. *Suppose that* $h$ *is concave and that the set* $\Lambda$ *is finite and nonempty. Then the semigroup of the solution of equation* (49) *preserves the function cone*

$$(86) \qquad \mathbf{F} = \left\{ f \in \mathbf{C}^2_{b1}(\mathbb{R}^\Lambda_{\geq 0}) \colon \frac{\partial}{\partial x_i} f \geq 0 \ \forall i, \frac{\partial^2}{\partial x_i \partial x_j} f \leq 0 \ \forall i, j \right\},$$

*where* $\mathbf{C}^2_{b1}(\mathbb{R}^\Lambda_{\geq 0})$ *denotes the space of all bounded* $\mathbf{C}^2$*- functions* $f \colon \mathbb{R}^\Lambda_{\geq 0} \to \mathbb{R}$ *with bounded first partial derivatives.*

Proof. This lemma is an addendum to Proposition 17 in [3]. There, the preservation of $\mathbf{F}$ was proved for $h \equiv 0$ and matrices $m$ with $\sum_{j \in \Lambda} m(i, j) = 1$ for all $i \in \Lambda$. This proof also works for more general matrices $m$ which only satisfy $\sum_{j \in \Lambda} m(i, j) \leq 1$ for $i \in \Lambda$. To extend the argument to the case $h \neq 0$, let $y(t, x)$ be the solution of

$$(87) \qquad \frac{\partial}{\partial t} y(t, x) = h(y(t, x)), \qquad y(0, x) = x \in I.$$

This defines a deterministic Markov process whose semigroup is given by $S_t f(x) := f(y(t, x))$. Similarly as in [3], we only need to establish that this semigroup preserves $\mathbf{F}$ if $h$ is twice continuously differentiable. A little calculation shows that it is enough to prove that $y(t, x)$ is increasing and concave in $x$. To show concavity, notice that differentiating equation (87) results in

$$(88) \qquad \frac{\partial}{\partial t} \frac{\partial^2}{\partial x^2} y(t, x) = h''(y(t, x)) \cdot \left( \frac{\partial}{\partial x} y(t, x) \right)^2 + h'(y(t, x)) \cdot \frac{\partial^2}{\partial x^2} y(t, x).$$

For fixed $x$, write (88) as $z'_t = a_t + b_t \cdot z_t$ with $z_0 = 0$. The solution for this is

$$(89) \qquad z_t = \exp\left( \int_0^t b_s \, ds \right) \int_0^t \exp\left( -\int_0^s b_r \, dr \right) a_s \, ds.$$

Since $h(x)$ is concave, $a_t$ is negative, implying the claimed concavity. A similar, even simpler argument shows monotonicity. □

PROOF OF PROPOSITION 2.2. We make use of the approximation scheme $X^\Lambda$ defined in the proof of Proposition 2.1; recall that $X^\Lambda$ is the solution of (49). Since $X_t^\Lambda \uparrow X_t$, it suffices to show the inequality (15) with $X_t$ replaced by $X_t^\Lambda$, and $f$ depending only on the coordinates $x_i$ with $i \in \Lambda$.

Furthermore, we assume for the rest of the proof that $h$ and $g$ satisfy Assumption A4. The general case follows then by approximating $h$ and $g$ pointwise by functions $h_k$ and $g_k$ satisfying Assumption A4. See Lemma 19 of [3] for the details.

In addition, we may assume that $f \in \mathbf{F}$; otherwise approximate $f$ by functions in $\mathbf{F}$ and use dominated convergence. Denote by $S_t$ the strongly continuous semigroup of $X^\Lambda$ defined on $\mathbf{C}(I^\Lambda)$. When applied to $\varphi \in \mathbf{F}$, the generator of $X^\Lambda$ takes the form (see, e.g., Theorem 7.3.3 of [19])

$$(90) \quad \mathcal{G}\varphi(x) = \sum_{i \in \Lambda} \left[ \alpha \left( \sum_{j \in \mathbb{N}} m(i,j) x_j - x_i \right) \frac{\partial}{\partial x_i} + h(x_i) \frac{\partial}{\partial x_i} + g(x_i) \frac{\partial^2}{\partial x_i^2} \right] \varphi(x).$$

By Proposition 1.1.5(c) of [8], we know that

$$(91) \quad \frac{d}{dt} S_t f = \mathcal{G} S_t f.$$

Let $\overline{V}_t = (\overline{V}_t(i))_{i \in \Lambda}$ be a system of processes coupled through the initial state $\overline{V}_0(i) = X_0(i)$, $i \in \Lambda$, but following independent mean field dynamics:

$$df(\overline{V}_t) = \alpha \sum_{i \in \Lambda} \frac{\partial}{\partial x_i} f(\overline{V}_t)(\mathbf{E}^{\bar{\mu}} \overline{V}_t(i) - \overline{V}_t(i) + h(\overline{V}_t(i))) \, dt$$

$$(92)$$

$$+ \sum_{i \in \Lambda} \frac{\partial^2}{\partial x_i^2} f(\overline{V}_t) g(\overline{V}_t(i)) \, dt + \sum_{i \in \Lambda} \frac{\partial}{\partial x_i} f(\overline{V}_t) \sqrt{2g(\overline{V}_t(i))} \, dB_t(i),$$

where $B(i)_{i \in \Lambda}$ are independent Brownian motions. Write $\bar{\mu}_t := \mathcal{L}(\overline{V}_t)$; for brevity, we suppress in this notation the dependence on $\Lambda$. By equation (92), the evolution of $\bar{\mu}_t$ is given by

$$\frac{d}{dt} \bar{\mu}_t f = \alpha \sum_{i \in \Lambda} \left[ \mathbf{E}^{\bar{\mu}} \left[ (\mathbf{E}^{\bar{\mu}} \overline{V}_t(i) - \overline{V}_t(i)) \left( \frac{\partial}{\partial x_i} f \right) (\overline{V}_t) \right] \right.$$

$$(93)$$

$$\left. + \bar{\mu}_t \left[ h(x_i) \frac{\partial}{\partial x_i} f + g(x_i) \frac{\partial^2}{\partial x_i^2} f \right] \right].$$

Integration by parts yields

$$(94) \quad \int_0^t \left( \frac{d}{ds} \bar{\mu}_s \right) S_{t-s} f \, ds = [\bar{\mu}_s S_{t-s} f]_0^t - \int_0^t \bar{\mu}_s \frac{d}{ds} S_{t-s} f \, ds.$$

In view of (91), this reads as

$$(95) \qquad \bar{\mu}_t f - \bar{\mu} S_t f = \int_0^t \left( \frac{d}{ds} \bar{\mu}_s - \bar{\mu}_s \mathcal{G} \right) S_{t-s} f \, ds.$$

We will show that the integrand is nonnegative. From Lemma 6.1, we know that $\varphi := S_{t-s} f$ (for $0 \le s \le t$ fixed) is an element of **F**. By equations (90) and (93),

$$
\begin{aligned}
(96) \qquad & \left( \frac{d}{ds} \bar{\mu}_s - \bar{\mu}_s \mathcal{G} \right) \varphi \\
& = \alpha \sum_{i \in \Lambda} \mathbf{E}^{\bar{\mu}} \overline{V}_t(i) \mathbf{E}^{\bar{\mu}} \left[ \left( \frac{\partial}{\partial x_i} \varphi \right)(\overline{V}_t) \right] \\
& \quad - \alpha \sum_{i,j \in \Lambda} m(i,j) \mathbf{E}^{\bar{\mu}} \left[ \overline{V}_t(j) \left( \frac{\partial}{\partial x_i} \varphi \right)(\overline{V}_t) \right] \\
& \ge \alpha \sum_{i \in \Lambda} \mathbf{E}^{\bar{\mu}} \overline{V}_t(i) \mathbf{E}^{\bar{\mu}} \left[ \left( \frac{\partial}{\partial x_i} \varphi \right)(\overline{V}_t) \right] \\
& \quad - \alpha \sum_{i,j \in \Lambda} m(i,j) \mathbf{E}^{\bar{\mu}} \overline{V}_t(j) \mathbf{E}^{\bar{\mu}} \left[ \left( \frac{\partial}{\partial x_i} \varphi \right)(\overline{V}_t) \right].
\end{aligned}
$$

Note that under the assumptions on $\bar{\mu}$ we have $\mathcal{L}(\overline{V}(i)) = \mathcal{L}(\overline{V}(0))$ for all $i \in \Lambda$. The right-hand side of (96) is nonnegative because of $\sum_{j \in \Lambda} m(i,j) \le 1$. To see the inequality in (96), notice that $-\frac{\partial}{\partial x_i} \varphi$ is bounded and componentwise increasing by Lemma 6.1. The claimed inequality thus follows from the fact that $\mathcal{L}(\overline{V}_t)$ is associated, which we now prove. Independent real-valued random variables are associated (see page 78 of [16]), and $\bar{\mu}$ is associated by assumption. Hence,

$$
\begin{aligned}
(97) \qquad \mathbf{E}^{\bar{\mu}}[f(\overline{V}_t) g(\overline{V}_t)] & = \mathbf{E}^{\bar{\mu}}[\mathbf{E}[f(\overline{V}_t) g(\overline{V}_t) | \overline{V}_0]] \\
& \ge \mathbf{E}^{\bar{\mu}}[\mathbf{E}[f(\overline{V}_t) | \overline{V}_0] \mathbf{E}[g(\overline{V}_t) | \overline{V}_0]] \\
& \ge \mathbf{E}^{\bar{\mu}}[\mathbf{E}[f(\overline{V}_t) | \overline{V}_0]] \mathbf{E}^{\bar{\mu}}[\mathbf{E}[g(\overline{V}_t) | \overline{V}_0]] \\
& = \mathbf{E}^{\bar{\mu}}[f(\overline{V}_t)] \mathbf{E}^{\bar{\mu}}[g(\overline{V}_t)],
\end{aligned}
$$

showing that $\mathcal{L}(\overline{V}_t)$ is associated. $\square$

PROOF OF THEOREM 1. Let $y_0 := \max\{y \ge 0 : h(y) = 0\}$. Assume for the moment that $y_0 > 0$. The measure $\bar{\mu} := \mathcal{L}(X_1^{(\infty)})$ is associated, shift invariant and its first moment is finite by Theorem 2. Let $(V_t)$ be the solution of (69) with initial distribution $\mu := \mathcal{L}(X_1^{(\infty)}(0))$. Theorem 2(e), (a) and Proposition 2.2 imply

$$(98) \qquad \mathbf{E} e^{-\lambda X_{t+1}(i)} \ge \mathbf{E} e^{-\lambda X_{t+1}^{(\infty)}(i)} = \mathbf{E}^{\bar{\mu}} e^{-\lambda X_t^{(\infty)}(i)} \ge \mathbf{E}^{\mu} e^{-\lambda V_t}.$$

It follows from Lemma 5.2 that, under the stated assumptions, $\mathbf{E}e^{-\lambda V_t} \to 1$ for all $\lambda > 0$ as $t \to \infty$. This proves the assertion for the case $y_0 > 0$.

If $y_0 = 0$, then $\tilde{h}(x) := (h(x) - h(1)) \wedge 0$ satisfies $h \leq \tilde{h} \leq 0$ because $h$ is concave. Let $\widetilde{X}$ be the solution of (2) with $h$ replaced by $\tilde{h}$ and with the same family of Brownian motions. By the previous step, $\widetilde{X}$ suffers local extinction. Lemma 3.3 implies $X \leq \widetilde{X}$, which completes the proof. $\square$

## 7. Self-duality. Proof of Theorem 3.

PROOF OF THEOREM 3. In the rest of the paper we exploit the specific form of the dynamics (1) for the interacting Feller diffusions with logistic growth. Theorem 3 states that the process is "dual to itself" via

$$(99) \qquad \mathbf{E}^x \exp\left(-\frac{\gamma}{\beta}\langle X_t, y\rangle\right) = \mathbf{E}^y \exp\left(-\frac{\gamma}{\beta}\langle x, X_t^\dagger\rangle\right).$$

We will prove this for the solution $X^\Lambda$ of (49). By (51), we know that the process $X^\Lambda$ monotonically approximates $X$. Hence, the assertion follows by dominated convergence.

For the rest of the proof, we consider $X^\Lambda$. We write $X$ instead of $X^\Lambda$ and $x, y$ instead of $x^\Lambda, y^\Lambda$. The duality function is $H(x,y) = \exp(-\frac{\gamma}{\beta}\langle x, y\rangle)$. Recall the definition of $\mathbf{C}_{b1}^2(\mathbb{R}_{\geq 0}^\Lambda)$ from Section 6. Define the linear operator $\mathcal{G}^X : \mathbf{C}_{b1}^2(\mathbb{R}_{\geq 0}^\Lambda) \to \mathbf{C}(\mathbb{R}_{\geq 0}^\Lambda)$ by

$$(100) \qquad \begin{aligned} \mathcal{G}^X f(x) = \sum_{i \in \Lambda}\Bigg[ &\alpha\left(\sum_{j \in \Lambda} m^\Lambda(i,j)x_j - x_i\right)\frac{\partial}{\partial x_i}f \\ &+ \gamma x_i(K - x_i)\frac{\partial}{\partial x_i}f + \beta x_i\frac{\partial^2}{\partial x_i^2}f \Bigg]. \end{aligned}$$

By Itô's formula, the process $(X_t)_t$ is a solution of the martingale problem for $(\mathcal{G}^X, \mathbf{C}_{b1}^2(\mathbb{R}_{\geq 0}^\Lambda))$. In order to apply Theorem 4.4.11 of [8] (with the choice $\alpha, \beta = 0$), we will show that

$$(101) \qquad \mathcal{G}^X H(\cdot, y)(x) = \mathcal{G}^{X^\dagger} H(x, \cdot)(y) \qquad \forall x, y \in \mathbb{R}_{\geq 0}^\Lambda.$$

We prove equation (101) by considering the different parts of (100) separately. Since $H$ is a function $\rho(\langle x, y\rangle)$ of the scalar product, it is easy to see that the migration terms of both sides are equal. To establish equation (101), it remains to show that

$$(102) \qquad \begin{aligned} &\gamma x_i(K - x_i)\frac{\partial}{\partial x_i}H(x,y) + \beta x_i\frac{\partial^2}{\partial x_i^2}H(x,y) \\ &= \gamma y_i(K - y_i)\frac{\partial}{\partial y_i}H(x,y) + \beta y_i\frac{\partial^2}{\partial y_i^2}H(x,y) \end{aligned}$$

for all $i \in \Lambda$. Observe that this equation is symmetric in $x$ and $y$. Consider the left-hand side of equation (102) divided by $H(x, y)$:

$$
\begin{aligned}
\gamma x_i(K - x_i) &\cdot \left(-\frac{\gamma}{\beta} y_i\right) + \beta x_i \cdot y_i^2 \left(\frac{\gamma}{\beta}\right)^2 \\
&= -\frac{\gamma^2 K}{\beta} x_i y_i + \frac{\gamma^2}{\beta} x_i^2 y_i + \frac{\gamma^2}{\beta} x_i y_i^2.
\end{aligned}
\tag{103}
$$

The right-hand side of (103) is symmetric in $x$ and $y$ and therefore, by interchanging the roles of $x$ and $y$, equals also the right-hand side of equation (102) divided by $H(x, y)$.

Theorem 4.4.11 of [8] is applicable if we prove that

$$
\sup_{s,t \leq T} |\mathcal{G}^X H(X_s, X_t^\dagger)|
\tag{104}
$$

is integrable for all $T < \infty$, where $X$ and $X^\dagger$ are independent. It is not hard to see that

$$
|\mathcal{G}^X H(x, y)| \leq C(|x||y| + |x| + |y|) \qquad \forall x, y \in \mathbb{R}_{\geq 0}^\Lambda
\tag{105}
$$

for a finite constant $C$. For this, use that $z \exp(-z)$ is bounded in $z \geq 0$. Integrability of (104) therefore follows from the independence of $X$ and $X^\dagger$ and from Lemma 3.2. $\quad \square$

Let us write $\mathcal{M}_c(\mathbb{Z}^d)$ for the set of configurations in $\mathbb{R}_{\geq 0}^{\mathbb{Z}^d}$ with finite support. As a consequence of the self-duality, we prove the following characterization of the upper invariant measure in terms of the finite mass process.

LEMMA 7.1. *Assume $\beta, \gamma > 0$. The upper invariant measure $\bar{\nu}$ of (1) is uniquely determined by*

$$
\int \exp\left(-\frac{\gamma}{\beta} \langle x, \lambda \rangle\right) \bar{\nu}(dx) = \mathbf{P}^\lambda(\exists t \geq 0 \text{ such that } X_t^\dagger = \underline{0}), \qquad |\lambda| < \infty,
\tag{106}
$$

*where $X^\dagger$ is the solution of (1) with the transpose migration matrix $m^\dagger$.*

PROOF. Fix a configuration $\lambda \in \mathcal{M}_c(\mathbb{Z}^d)$ and consider the process $(X_t^{(n)})$ started in the constant configuration $\underline{n}(i) \equiv n$. This process converges to the maximal process as $n \to \infty$. Therefore, the self-duality implies, for $t > 0$,

$$
\begin{aligned}
\mathbf{E} \exp\left(-\frac{\gamma}{\beta} \langle X_t^{(\infty)}, \lambda \rangle\right) &= \lim_{n \to \infty} \mathbf{E} \exp\left(-\frac{\gamma}{\beta} \langle X_t^{(n)}, \lambda \rangle\right) \\
&= \mathbf{E}^\lambda \lim_{n \to \infty} \exp\left(-\frac{\gamma}{\beta} \langle X_t^\dagger, \underline{n} \rangle\right) \\
&= \mathbf{P}^\lambda(X_t^\dagger = \underline{0}) \\
&= \mathbf{P}^\lambda(\exists s \leq t : X_s^\dagger = \underline{0}).
\end{aligned}
\tag{107}
$$

For the second equality, we used monotone convergence. Letting $t \to \infty$, the assertion follows from Theorem 2(d). For general $\lambda$ with $|\lambda| < \infty$, use monotone convergence. □

## 8. Convergence to the upper invariant measure. Proof of Theorem 5.

PROOF OF THEOREM 5. Let $\mu$ be a translation invariant distribution on $\mathbb{Z}^d$ with $\mu(\underline{0}) = 0$. For analyzing the long-term behavior of the interacting Feller diffusion with logistic growth started in $\mu$, we can assume without loss of generality that $\mu$ has finite first moment and satisfies $\mu(x_0 = 0) = 0$. Otherwise, we let the system run for a little time $\varepsilon > 0$, obtaining

$$\lim_{t \to \infty} \mathcal{L}^{\mu}(X_t) = \lim_{t \to \infty} \mathcal{L}^{\mathcal{L}^{\mu}(X_\varepsilon)}(X_t). \tag{108}$$

A comparison with the maximal process [see Theorem 2(e), (b)] yields $\mathbf{E}^{\mu} X_\varepsilon(0) < \infty$. Furthermore, after a fixed positive time $\varepsilon > 0$, every component is strictly positive almost surely (see Lemma 8.2).

Let $X$ and $X^{\dagger}$ be solutions of (1) with migration matrix $m$ and its transpose $m^{\dagger}$, respectively. In Lemma 8.1 we will show that the total mass hits zero in finite time or tends to infinity. Hence, we get by self-duality (Theorem 3)

$$
\begin{aligned}
\mathbf{E}^{\mu} &\exp\left(-\frac{\gamma}{\beta}\langle X_t, \lambda \rangle\right) \\
&= \int \mu(dx) \left( \mathbf{E}^{\lambda} \left[ 1_{|X_s^{\dagger}| \to \infty} \exp\left(-\frac{\gamma}{\beta}\langle x, X_t^{\dagger} \rangle\right) \right] \right. \\
&\qquad\qquad \left. + \mathbf{E}^{\lambda} \left[ 1_{\exists s : X_s^{\dagger} = \underline{0}} \exp\left(-\frac{\gamma}{\beta}\langle x, X_t^{\dagger} \rangle\right) \right] \right).
\end{aligned}
\tag{109}
$$

We treat the two terms on the right-hand side separately and begin with the first term. Apply Hölder's inequality to the integral with respect to $\mu$. For this, let $1/p_i = X_t(i)/|X_t|$ if this is positive. Thus, we obtain

$$
\begin{aligned}
\mathbf{E}^{\lambda} &\left[ 1_{|X_s^{\dagger}| \to \infty} \int \mu(dx) \exp\left(-\frac{\gamma}{\beta}\langle x, X_t^{\dagger} \rangle\right) \right] \\
&\leq \mathbf{E}^{\lambda} \left[ 1_{|X_s^{\dagger}| \to \infty} \prod_{i \in \mathbb{Z}^d} \left( \int \mu(dx) \exp\left(-\frac{\gamma}{\beta} x_i |X_t^{\dagger}|\right) \right)^{X_t(i)/|X_t|} \right] \\
&= \mathbf{E}^{\lambda} \left[ 1_{|X_s^{\dagger}| \to \infty} \int \mu(dx) \exp\left(-\frac{\gamma}{\beta} x_0 |X_t^{\dagger}|\right) \right] \to 0 \qquad \text{as } t \to \infty.
\end{aligned}
\tag{110}
$$

The equality is a consequence of the translation invariance of $\mu$. The last expression tends to zero because of dominated convergence and the assumption $\mu(x_0 = 0) = 0$. As to the second term on the right-hand side of (109), dominated convergence gives

$$(111) \quad \int \mu(dx) \mathbf{E}^\lambda \left[ 1_{\exists s \, : \, X_s^\dagger = \underline{0}} \exp\left( -\frac{\gamma}{\beta} \langle x, X_t^\dagger \rangle \right) \right] \to \int \mu(dx) \mathbf{E}^\lambda [1_{\exists s \, : \, X_s^\dagger = \underline{0}}]$$

as $t \to \infty$. Using Lemma 7.1, we arrive at

$$\lim_{t \to \infty} \mathbf{E}^\mu \exp\left( -\frac{\gamma}{\beta} \langle X_t, \lambda \rangle \right) = \mathbf{P}^\lambda (\exists t \geq 0 \text{ s.t. } X_t^\dagger = \underline{0}) = \int \exp\left( -\frac{\gamma}{\beta} \langle x, \lambda \rangle \right) \bar{\nu}(dx).$$

Starting in $\mathcal{L}(X_0) \geq \mu$, the process $\mathcal{L}(X_t)$ is bounded below by $\mathcal{L}^\mu(X_t)$ (Lemma 3.3) and is bounded above by $\mathcal{L}(X_t^{(\infty)})$ [Theorem 2(e)], which both converge to $\bar{\nu}$. This concludes the proof of Theorem 5. $\square$

We have to append the following:

LEMMA 8.1. *Assume* $\beta > 0$. *Let* $X$ *be a solution of* (1) *starting in* $x \in \mathbb{E}_\sigma$ *with finite total mass* $|x| < \infty$. *Then with probability* 1, *either,*

- *there is a* $t \geq 0$ *such that* $X_s = \underline{0}$ *for all* $s \geq t$ *or*
- $|X_t| \to \infty$ *as* $t \to \infty$.

PROOF. The intuition behind this is the following. The process always has a positive probability of hitting the lower trap. Whenever the total mass stays bounded, the process will seize its chance.

This is made precise in Theorem 2 of [13]. In order to apply this result, we only need to verify that there always is the risk of extinction in the following sense:

$$(112) \quad \forall y \colon \inf_{|x| \leq y} \mathbf{P}^x (\exists t \colon X_t = \underline{0}) > 0.$$

Let $Y_t$ be a solution of (2) with $G = \mathbb{Z}^d$, $h(x) = \gamma K x$ and $g(x) = \beta x$. By Lemma (3.3), $X$ and $Y$ may be coupled such that $X_t$ is bounded above by $Y_t$ almost surely. Furthermore, $|Y_t|$ is equal in distribution to Feller's branching diffusion $F_t$ with super-criticality $\gamma K$ started in $|x|$. The extinction probability of $F_t$ is strictly positive; see, for example, Appendix 6.2 of [6]. Therefore, condition (112) follows from

$$(113) \quad \mathbf{P}^x (\exists t \geq 0 \colon X_t = \underline{0}) \geq \mathbf{P}^{|x|} (\exists t \geq 0 \colon F_t = 0) \geq \mathbf{P}^y (\exists t \geq 0 \colon F_t = 0) > 0$$

for every $x$ with $|x| \leq y$. $\square$

LEMMA 8.2. *Suppose that* $h$ *and* $g$ *satisfy Assumption* A1. *Let* $X$ *be a solution of* (2). *If its initial law* $\mu$ *is translation invariant and does not charge the zero configuration* $\underline{0}$, *then, for every fixed time* $t_0 > 0$,

$$(114) \quad X_{t_0}(i) > 0 \qquad \forall i \in G \ \mathbf{P}^\mu\text{-}a.s.$$

PROOF. Assume that $h \leq 0$. Otherwise, compare $X_t$ with the process defined with $h \wedge 0$ instead of $h$.

Let $\tilde{h}(\theta) = -\alpha\theta + h(\theta)$. For $\varepsilon > 0$, define the solution of

$$(115) \qquad dY_t^{\varepsilon,i} = \alpha\varepsilon\,dt + \tilde{h}(Y_t^{\varepsilon,i})\,dt + \sqrt{2g(Y_t^{\varepsilon,i})}\,dB_t^i, \qquad Y_0^{\varepsilon,i} \geq 0,$$

on the same probability space as $X$ by using the same system of Brownian motions. This system satisfies $\mathbf{P}^0(Y_t^{\varepsilon,i} > 0) = 1$ for all $t > 0$. Otherwise, continuity in the initial value would imply that there is a $t > 0$ and a $\theta_0$ such that $\mathbf{P}^\theta(Y_t^{\varepsilon,i} = 0) > 0$ for all $\theta \leq \theta_0$. Integrating this with the equilibrium distribution $\Gamma_\varepsilon$ [see equation (74), it exists because of $h \leq 0$] yields $\Gamma_\varepsilon(0) > 0$, which is false. Thus, we have

$$(116) \qquad \begin{aligned} &Y_\delta^{\varepsilon m(i,j),i} > 0 \\ &\qquad \forall \varepsilon \in (0,1) \cap \mathbb{Q} \qquad \forall \delta \in (0,1) \cap \mathbb{Q} \qquad \forall i,j \text{ s.t. } m(i,j) > 0 \qquad \text{a.s.} \end{aligned}$$

Denote the event $\{X_t(j) \geq \varepsilon \forall t \in [t_0 - \delta, t_0]\}$ by $A_{\varepsilon,\delta}$. On $A_{\varepsilon,\delta}$ we compare $X$ with the solution of (115):

$$(117) \qquad \begin{aligned} X_t(i) &= X_{t_0-\delta}(i) + \int_{t_0-\delta}^t \alpha \sum_k m(i,k)X_s(k)\,ds \\ &\quad + \int_{t_0-\delta}^t \tilde{h}(X_s(i))\,ds + \int_{t_0-\delta}^t \sqrt{2g(X_s(i))}\,dB_s(i) \\ &\geq \int_{t_0-\delta}^t \alpha m(i,j)\varepsilon\,ds + \int_{t_0-\delta}^t \tilde{h}(X_s(i))\,ds \\ &\quad + \int_{t_0-\delta}^t \sqrt{2g(X_s(i))}\,dB_s(i), \end{aligned}$$

for all $t \in [t_0 - \delta, t_0]$. By standard comparison results (e.g., Theorem (V.43.1) in [20] and a stopping argument), this implies $X_{t_0}(i) \geq Y_\delta^{m(i,j)\varepsilon}$ on $A_{\varepsilon,\delta}$ a.s. By path continuity, $A_{\varepsilon,\delta}$ approximates $\{X_{t_0}(j) > 0\}$ as $\delta, \varepsilon \to 0$. It follows that on $X_{t_0}(j) > 0$ we have $X_{t_0}(i) > 0$ for all $i$ such that $m(i,j) > 0$ a.s. With the migration kernel being irreducible, every site can be reached from $j$. By induction, we conclude that every component of $X_{t_0}$ is positive a.s. given $X_{t_0}(j) > 0$.

Starting in a nontrivial translation invariant measure, the system a.s. never hits $\underline{0}$. Therefore, there is a location $j$ with $X_{t_0}(j) > 0$ a.s. This proves the lemma. □

**Acknowledgments.** We thank Frank den Hollander, Jan Swart and, in particular, Don Dawson for valuable discussions. Also, we thank two referees for a number of very helpful suggestions.

INSTITUT FÜR STOCHASTIK
    UND MATHEMATISCHE INFORMATIK
UNIVERSITÄT FRANKFURT
ROBERT-MAYER-STR. 6-10
60325 FRANKFURT/MAIN
GERMANY
E-MAIL: hutzenth@math.uni-frankfurt.de
        wakolbin@math.uni-frankfurt.de